\documentclass[11pt]{amsart}
\usepackage{amssymb, amscd}
\usepackage{latexsym,epsfig}
\usepackage[all]{xy}
\usepackage{pst-all}
\numberwithin{equation}{section}

\newcommand{\VVl}{\VV_{\lambda}}

\newcommand{\VV}{\mathbb V}
\newcommand{\LL}{\mathbb L}
\newcommand{\EE}{\mathbb E}

\newcommand{\PP}{\mathbb P}
\newcommand{\DD}{\mathbb D}
\newcommand{\GG}{\mathbb G}

 \newcommand{\CC}{\mathbb C}

\newcommand{\QQ}{\mathbb Q}
\newcommand{\ZZ}{\mathbb Z}

\newcommand{\G}{\bf G}

\newcommand{\Ag}{{\mathcal A}_g}

\newcommand{\Agp}{{\mathcal A}_g^{\prime}}
\newcommand{\Agf}{\tilde{\mathcal A}_g}
\newcommand{\Xg}{{\mathcal X}_g}
\newcommand{\Xgmin}{{\mathcal X}_{g-1}}

\newtheorem{theorem}{Theorem}[section]
\newtheorem{lemma}[theorem]{Lemma}
\newtheorem{proposition}[theorem]{Proposition}
\newtheorem{corollary}[theorem]{Corollary}

\newtheorem{definition-lemma}[theorem]{Definition-Lemma}

\theoremstyle{definition}

\newtheorem{example}[theorem]{Example}

\theoremstyle{remark}
\newtheorem{remark}[theorem]{Remark}

\begin{document}

\title[]{Rank One Eisenstein Cohomology 
\\ of Local Systems on the moduli space
\\ of abelian varieties}
\author{Gerard van der Geer}
\address{Korteweg-de Vries Instituut, Universiteit van
Amsterdam, Plantage 
\newline Muidergracht 24, 1018 TV Amsterdam, The Netherlands.}
\email{geer@science.uva.nl}


\begin{abstract}
We give a formula for the Eisenstein cohomology of local systems
on  the partial compactification of the moduli of principally polarized 
abelian varieties given by rank~$1$ degenerations.
\end{abstract}

\maketitle

\begin{section}{Introduction}
\label{sec: intro}
Let $\Ag$ be the moduli stack of principally polarized abelian varieties
over a field $k$. 
The universal abelian scheme $\pi: \Xg \to \Ag$ over $\Ag$ defines a local 
system $\VV= R^1\pi_*\QQ(1)$ of rank $2g$ on $\Ag$, and 
for a prime $\ell$ different from the characteristic of $k$
its $\ell$-adic variant
$\VV= R^1\pi_*\QQ_{\ell}(1)$ for the \'etale topology. 
This local system is associated to the standard representation of 
${\rm GSp}(2g,\QQ)$. To an
irreducible representation of ${\rm Sp}(2g,\QQ)$ with highest weight $\lambda
=(\lambda_1,\ldots,\lambda_g)$ with $\lambda_1\geq \lambda_2 \geq \ldots
\geq \lambda_g$ we can associate a local system $\VVl$  of weight
$\sum_i \lambda_i$ which occurs `for the first time' in
$$\otimes_{j=1}^{g-1} {\rm Sym}^{\lambda_j-\lambda_{j+1}}(\wedge^j \VV) 
\otimes  \wedge^g \VV^{\lambda_g}.
$$
The cohomology of
this local system is closely related to vector-valued Siegel modular
forms, cf.\ \cite{Del,F-C,FvdG}. 
We are interested in the Euler characteristic
$$
e(\Ag,\VVl):=\sum (-1)^i \, [H^i(\Ag,\VVl)]
$$
in a suitable $K$-group of mixed Hodge structures or
of Galois representations if one works with \'etale cohomology. 
Similarly, we can consider
the analogue for compactly supported cohomology
$$
e_c(\Ag,\VVl):=\sum (-1)^i \, [H_c^i(\Ag,\VVl)].
$$
There is a natural map $H_c^*(\Ag,\VVl) \to H^*(\Ag,\VVl)$ and the image is
called the interior cohomology. We define the Eisenstein cohomology as the
difference
$$
e_{\rm Eis}(\Ag,\VVl):=e(\Ag,\VVl)-e_c(\Ag,\VVl).
$$
We consider 
the partial compactification of rank $\leq 1$ degenerations. One can obtain it
by considering a compactification $\tilde{\Ag}$ of Faltings-Chai-type
together with the natural
map $q: \tilde{\Ag} \to {\mathcal A}_g^*$ to the Satake 
compactification. The Satake compactification has a stratification
${\mathcal A}_g^* =\cup_{i=0}^g{\mathcal A}_{i}$. 
Then the partial compactification ${\Agp}$ is defined
as $q^{-1}(\Ag \cup {\mathcal A}_{g-1})$ and can be  viewed as
the moduli space of semi-abelian varieties with torus part of rank $\leq 1$.
It is independent of the choice of the compactification $\tilde{\Ag}$.

Let $j: {\Ag} \to {\Agf}$ be the inclusion map. 
We can rewrite the Eisenstein cohomology for a local system $\VVl$
(both in Hodge cohomology and in \'etale cohomology)
as 
$$
e_{\rm Eis}(\Ag,\VVl)=e(\Agf,Rj_*\VVl)-e(\Agf,Rj_{!}\VVl),
$$ 
where $j_!\VVl$ is the extension by zero. This is a sum over strata 
$$
\sum_{j=0}^{g-1} e_c(q^{-1}({\mathcal A}_{j}), Rj_*\VVl-Rj_!\VVl).
$$

We define now the
rank $1$ part of the Eisenstein cohomology as the part that comes from
the boundary component lying over ${\mathcal A}_{g-1}$:
$$
e_{\rm Eis,1}({\mathcal A}_g, \VVl)
:=e_c(q^{-1}({\mathcal A}_{g-1}), Rj_*\VVl-Rj_!\VVl).
$$
This contribution to the Eisenstein cohomology is independent of the
compactification.

Let $\LL:=h^2({\PP}^1)$ be the Lefschetz motive $h^2({\PP}^1)$ of rank $1$ 
and weight~$2$. Our result is an explicit formula for the rank $1$ part
of the Eisenstein cohomology.

\begin{theorem}\label{rank1thm}
The contribution $e_{\rm Eis,1}(\Ag,\VVl)$
to the Eisenstein cohomology of $\VVl$
from the codimension $1$ boundary 
is (both in Hodge cohomology and in \'etale cohomology)
of the form
$$
\sum_{k=1}^g (-1)^{k} \, e_c({{\mathcal A}_{g-1}},
\VV_{\lambda_1+1,\lambda_2+1,\ldots,\lambda_{k-1}+1,
\lambda_{k+1}, \ldots,\lambda_g})\,\, (1-\LL^{\lambda_{k}+g+1-k}).
$$
\end{theorem}
\begin{remark}
The fact that the result is formally the same for both Hodge and
\'etale cohomology can be deduced by using the compactifications
of the powers of the universal abelian variety as in \cite{F-C}
and by decomposing the direct image of the cohomology under the
action of algebraic correspondences as done in \cite{F-C}, p.\ 238.
\end{remark}
\begin{example}

For $g=1$ the cohomology of $\VV_k$ can only be non-trivial for even $k$.
In this case on gets for the Eisenstein cohomology 
$e_{\rm Eis}({\mathcal A}_1,\VV_k)=e_{\rm Eis,1}({\mathcal A}_1,\VV_k)$
of 
$\VV_{k}$ the polynomial $1-\LL^{l+1}$,
in agreement with the Eichler-Shimura isomorphisms of \cite{Del}
$$
e_c({\mathcal A}_1,\VV_{k})
= -S_{k+2} \oplus \bar{S}_{k+2}-1
$$
and
$$
e({\mathcal A}_1,\VV_{k})=-S_{k+2} \oplus \bar{S}_{k+2}-\CC(k+1).
$$
Here $S_k$ denotes the space of cusp forms of weight $k$ on ${\rm SL}(2,\ZZ)$.

For $g=2$, cohomology of $\VV_{l,m}$ can only be non-trivial for 
$l\equiv m \, (\bmod 2)$ and in this case one gets for
$e_{{\rm Eis,1}}(\VV_{l,m})$ the expression
$$
e_c({\mathcal A}_1,\VV_{l+1}) (1-\LL^{m+1}) 
- e_c({\mathcal A}_1,\VV_m)(1-\LL^{l+2}),
$$
and for $g=3$ we get for $e_{{\rm Eis},1}$ the expression
\begin{multline*}
e_c({\mathcal A}_2, \VV_{l+1,m+1})(1-\LL^{n+1})
-e_c({\mathcal A}_2,\VV_{l+1,n})(1-\LL^{m+2})\\
+e_c({\mathcal A}_2,\VV_{m,n})(1-\LL^{l+3}).
\end{multline*}
For $g=2$ we also prove a formula for the {\sl total} Eisenstein cohomology
$$
\begin{aligned}
e_{\rm Eis}({\mathcal A}_2,\VV_{l,m})&=
-s_{l-m+2}(1-\LL^{l+m+3})+  s_{l+m+4}(\LL^{m+1}-\LL^{l+2})+ \qquad \\
&\qquad +  \begin{cases}
e_c({\mathcal A}_1,\VV_m)(1-\LL^{l+2})-(\LL^{l+2}-\LL^{l+m+3})   
& \text{ $l$ even,}\cr
-e_c({\mathcal A}_1,\VV_{l+1})(1-\LL^{m+1})-(1-\LL^{m+1}) & \text{$l$ odd,}\cr
\end{cases} \\
\end{aligned}
$$
where $s_m$ denotes for $m>2$ the dimension of the space of cusp 
forms of weight $m$ 
on ${\rm SL}(2,\ZZ)$ and $s_2=-1$. 
From this formula one can deduce for regular $\lambda$
(i.e., $l>m>0$) the formula for the Eisenstein cohomology that was
announced in joint work with Carel Faber \cite{FvdG}, see Corollary
\ref{gistwee}. In that paper the term Eisenstein cohomology refers only to the 
kernel of $H_c^* \to H^*$.
\end{example}

\begin{remark} Note that the result of \ref{rank1thm} is compatible
with Poincar\'e duality,  which says that
$$
H^i(\Ag,\VVl)^{\vee} \cong H^{2d-i}_c(\Ag,\VVl^{\vee}(\nu^{d})),
$$
where $d=g(g+1)/2$. Here ${\VV}(\nu)$ means the twist of $\VV$
by the multiplier, cf.\ Section \ref{TheGroup}.
\end{remark}
The study of Eisenstein cohomology was initiated by Harder and 
carried on by his students Schwermer and Pink, cf.\
\cite{HarderSLN, Schwermer1,Schwermer2, Pink},
cf.\ also the work of Franke \cite{Franke}. One may view the result
here as an explicit formula for the general results of 
\cite{Pink} for the symplectic group. 
For us the interest in  Eisenstein cohomology
arose in joint work \cite{FvdG} with Carel Faber 
where we tried to obtain information on Siegel modular
forms by counting curves over finite fields. There Eisenstein cohomology
contributes terms that one wants to remove, cf.\ \cite{BFG,FvdG}.

The author thanks Gerd Faltings for some useful remarks.

\end{section}
%%%
%%%
\begin{section}{The Group}\label{TheGroup}
Fix a positive integer $g$ and let $V_{\ZZ}$ be the standard
symplectic lattice of rank $2g$ with generators $e_i$ and $f_i$
($i=1,\ldots,g$) with $\langle e_i, f_j\rangle = \delta_{ij}$,
the Kronecker delta,
for $i \leq j$. We let $G={\rm GSp}(2g)$ be the corresponding
Chevalley group of symplectic similitudes of $V_{\ZZ}$.
We shall write $V$ for ${\VV}_{\ZZ}\otimes {\QQ}$.
 
An element $\gamma \in G$ can be written as 
$
\gamma =
\left(\begin{matrix} a & b\cr
c & d\cr\end{matrix}\right) \in G
$
with $a,b,c,d$ integral $g\times g$-matrices. For such a $\gamma$
we write $a\, d^t - b \, c^t= \nu(\gamma) \, {\rm Id}_g$.
Here $\nu : G \to {\G}_m$ is called the multiplier representation.
It satisfies $\det = \nu^g: G \to {\G}_m$. We let $G$ act on the left
on $V$ by matrix multiplication.
We denote by $M$ the subgroup of $G$ of elements that respect
the two subspaces $\langle e_i: i=1,\ldots g\rangle$ and
$\langle f_i: i=1,\ldots g\rangle$. 
We can interpret the
elements of $M$ as matrices $\left(\begin{matrix} a & 0\cr
0 & d\cr\end{matrix}\right)$ with $a\cdot d^t = \nu \, 1_g$.
So $M \cong {\rm GL}(g) \times {\G}_m$.
 Furthermore, we let $Q$ be the
(maximal) parabolic
subgroup of $G$ that stabilizes the sublattice spanned by
the vectors $e_i$ ($i=1,\ldots,g$). Then $Q=M \ltimes U$
with $U={\rm Hom}({\rm Sym}^2(X), \ZZ)$, the group of $\ZZ$-valued
bilinear forms, or in terms of matrices, the group of matrices
$\left(\begin{matrix} 1_g & b\cr
0 & 1_g\cr\end{matrix}\right)$ in $G$.
The standard maximal torus $T$ of $G$ can be identified with
${\GG}_m^{g+1}$ via
$$
(t_1,\ldots,t_g,x) \mapsto {\rm diag}(t_1,\ldots,t_g,x/t_1,\ldots,x/t_g).
$$
The character group may thus be identified with the lattice
$$
\{ (a_1,\ldots,a_g,c) \in {\ZZ}^{g+1} : \sum a_i \equiv c \, (\bmod 2)\},
$$
via $\chi({\rm diag}(t_1,\ldots,t_g,x/t_1,\ldots,x/t_g))= 
x^{(c-\sum a_i)/2}\prod t_i^{a_i}$.
We will sometimes view $t_i$ ($i=1,\ldots,g$) and $x$ as characters on $T$.

Let $t$ be the complex Lie algebra of $T$. We have root
systems $\Phi_G$ and $\Phi_M$ in $t^{\vee}$ and we choose
compatible systems of positive roots $\Phi_G^{+}$ and $\Phi_M^{+}$.
So
$$
\Phi_M:=\{ (t_i/t_j)^{\pm} : 1 \leq i < j \leq g \},
$$
and we set
$$
\Psi_M= \{ (t_it_j/x)^{\pm} : 1 \leq i \leq j \leq g \}
$$
so that
$$
\Phi_G=\Phi_M \cup \Psi_M \qquad \hbox{\rm and} \quad
\Phi_G^+=\Phi_M^+ \cup \Psi_M^+
$$
with
$$
\Phi_M^{+}=\{ (t_i/t_j) : 1 \leq i < j \leq g \}
\quad \hbox{\rm and} \quad
\Psi_M^{+}= \{ (t_it_j/x) : 1 \leq i \leq j \leq g \}.
$$
As usual $\rho=(1/2)\sum_{x \in \Phi_G^{+}} x$.
\end{section}
\begin{section}{The Final Elements of the Weyl Group}\label{Final}
The Weyl group $W_G$ of $G$ is isomorphic to the semi-direct product
$S_g \ltimes ({\ZZ} / 2{\ZZ})^g$, where the symmetric group
$S_g$ on $g$ letters acts on $({\ZZ} / 2{\ZZ})^g$ by permuting the $g$
factors. We interpret elements of $W_G$ as signed permutations.
The Weyl group $W_M$ of $M$ is isomorphic to the symmetric
group $S_g$. They operate on the complex Lie algebra $t$ and 
its dual $t^{\vee}$.

We define the set of Kostant representatives by
$$
W^M= \{ w \in W_G : \Phi_M^{+} \subset w(\Phi_G)^{+} \},
$$
or equivalently as
$$
\{ w \in W_G : \langle w(\rho)-\rho, u\rangle \geq 0
\quad \hbox{\rm for all $u \in \Phi_M^{+}$} \}.
$$
Concretely, for a sign change $\epsilon: (x_1,\ldots,x_g)\mapsto
(\epsilon_1x_1,\ldots,\epsilon_gx_g)$ with $\epsilon_i \in \{\pm 1\}$
there exists exactly one $\sigma \in S_g$ such that
$\sigma\epsilon (\rho)-\rho$ is of the form $(a_1,\ldots,a_g)$
with $a_1\geq a_2 \geq \cdots \geq a_g$ and this is the Kostant
representative in $W^M$.
Note that $\rho=(g,g-1,\ldots, 2,1,0)$. Recall that $W_G$ carries a
length function $\ell$.

Another way to describe the Weyl group (and the one we shall use in 
the following) is as the group  of permutations
$$
W_g:= \{ \sigma \in S_{2g} : \sigma(i)+\sigma(2g+1-i)=2g+1 \, (i=1,\ldots,g)\}.
$$
Then the Weyl group of $M$ can be identified with the subgroup
$$
S_g=\{\sigma \in W_g:
\sigma\{1,2,\ldots,g\} =\{1,2,\ldots,g\}\}.
$$
The length $\ell(w)$ of an element $w\in W_g\subset S_{2g}$ 
is then defined by
$$
\#\{i < j \le g : w(i) > w(j) \} + \# \{ i \leq j \le g : w(i)+w(j)>
2g+1\}.
$$
\begin{lemma}
An element $\sigma \in W_g\subset S_{2g} $ is a Kostant representative 
if and only if $\sigma(i)<\sigma(j)$ for all $1\leq i < j \leq g$.
\end{lemma}
For the proof we refer to  \cite{EvdG}, Lemma 1.
In that paper the Kostant representatives are called {\sl final elements}.
We shall adopt that usage here too.
The set of $2^g$ final elements in the Weyl group
is denoted by $F_g\subset S_{2g}$. 

\begin{lemma}\label{subperms}
 The set $\{\sigma \in F_g : \sigma^{-1}(k)\leq k \}$ (resp.\
the set $\{\sigma \in F_g : \sigma^{-1}(2g+1-k)\leq g \}$)
has cardinality $2^{g-1}$ and can be identified in a natural way
with $F_{g-1}$ compatibly with the length function $\ell$.
\end{lemma}
\begin{proof}
Identify an element $\sigma$ with its image $g$-tuple
 $[\sigma(1),\ldots,\sigma(g)]$.
Let $A = \{ \sigma \in F_g : \sigma^{-1}(k) \leq k\}$. For each $\sigma \in A$
we delete the entry equal to $k$ from
$[\sigma(1),\ldots,\sigma(g)]$. We rename the entries by replacing an
entry $m$ by $m-1$ if $k<m < 2g+1-k$ and by $m-2$ if $m>2g+1-k$. We thus find
the $2^{g-1}$ elements of $F_{g-1}$ as one easily checks. The other
statements are proved in a similar way.
\end{proof}
\begin{example}
{\bf $g=3$} (For the third column see next two sections.)
\begin{center}
\begin{tabular}[t]{|c|c|c|c|}
\hline
$\ell(w)$ & $[w(1),w(2),w(3)]$ & $w(\lambda+\rho)-\rho$ \\
\hline
$0$ & $[123]$ & $(l,m,n)$ \\
$1$ & $[124]$ & $(l,m,-n-2)$ \\
$2$ & $[135]$ & $(l,n-1,-m-3)$ \\
$3$ & $[236]$ & $(m-1,n-1,-l-4)$ \\
$3$ & $[145]$ & $(l,-n-3,-m-3)$ \\
$4$ & $[246]$ & $(m-1,-n-3,-l-4)$ \\
$5$ & $[356]$ & $(n-2,-m-4,-l-4)$ \\
$6$ & $[456]$ & $(-n-4,-m-4,-l-4)$ \\
\hline
\end{tabular}
\end{center}
\end{example}
\end{section}
\begin{section}{Representations and Vector Bundles}\label{Reps}
Let $t$ be the complex Lie algebra of $T$.
The irreducible representations of $G$ are parametrized by the
characters of $T$ that correspond to $G$-dominant (i.e.,
the scalar product with $\Phi_G^{+}$ is non-negative) integral
weights $\lambda \in t^{\vee}$.  If $\lambda$
is given by $(a_1,\ldots,a_g,c)$ this means that we have
$a_1\geq a_2 \geq \ldots \geq a_g$. The standard representation
corresponds to $(1,0,\ldots,0,1)$. The irreducible representation
associated to $\lambda$ is denoted by $V(\lambda)$.

We have
$$
V(\lambda)^{\vee}= V(\lambda)\otimes \nu^k
\text{ with $k$ equal to $\lambda$ evaluated at $-1_{2g} \in t$}.
$$
The irreducible representations of $M$ are parametrized by characters $\mu$
of $T$ which are $M$-dominant (i.e., with non-negative scalar
product with $\Phi_M^{+}$) so that the representation $W(\mu)$ corresponding
to $\mu$ has highest weight $\mu$.

\smallskip
Let $\Ag$ be the Deligne-Mumford stack of principally
polarized abelian varieties of dimension $g$ and let $\pi :
\Xg \to \Ag$ be the universal family.
In the complex category we may identify $\Ag({\CC})$
with ${\rm Sp}(2g,{\ZZ})\backslash {\DD}$ with ${\DD}$
the space of Lagrangian subspaces $W \subset V_{\CC}$
on which $-\langle v,\bar{v}\rangle >0$. The complex
manifold ${\DD}$ is contained in its so-called
compact dual ${\DD}^{\vee}= G/Q({\CC})$ of Lagrangian
subspaces $W \subset V_{\CC}$.

To each finite-dimensional complex representation $r$ of $Q({\CC})$
on a vector space $U$ one can associate a 
$G$-equivariant bundle $E^{\prime}_r$
on ${\DD}^{\vee}$ defined as the quotient of $G({\CC})\times U$
under the equivalence relation
$$
(gq,r(q)^{-1} u) \sim (g,u) \qquad
\hbox{\rm for all $g \in G(\CC)$, $q \in Q(\CC)$}.
$$
Its restriction to ${\DD}$
descends to a vector bundle $E_r$ on $\Ag({\CC})$.
If $r$ is the restriction to $Q({\CC})$ of a finite-dimensional
complex representation $\rho$ of $G({\CC})$ then $E_r$ carries
a integrable connection defined as follows. An element $u \in U$,
the fibre over the base point, defines a trivialization
of $E_r$ by $\gamma \mapsto \rho(\gamma)u$; hence an
integrable connection on $E_r$ and it descends to $\Ag({\CC})$.

The vector bundle associated to the standard representation of $G(\CC)$
can be identified with the relative de Rham homology of
$\Xg$ and $R^1\pi_*({\CC})$ is the vector bundle
associated to $\nu^{-1} \otimes$ the standard representation.
The integrable connection is the Gauss-Manin connection.

After the choice of a base point we can view $G$ as the
fundamental group (arithmetic fundamental group) of $\Ag$.
Therefore we can associate a local system (or a smooth ${\QQ}_l$-sheaf)
to each finite-dimensional representation of $G$.
To the standard representation it associates the local
system $R^1\pi_*{\QQ}(1)$ (resp.\ ${\QQ}_{l}$-sheaf
$R^1\pi_*{\QQ}_l(1)$ on $\Xg \otimes {\ZZ}[1/l]$).
The character $\nu$ corresponds to ${\QQ}(1)$ (or ${\QQ}_{\ell}(1)$).
We have a non-degenerate alternating pairing
$$
R^1\pi_*{\QQ} \times R^1\pi_*{\QQ} \to {\QQ}(-1)
$$

The local system associated
to such a $\lambda$ is denoted by $\VVl$. If $\lambda$
is given by $(a_1,\ldots,a_g,c)$ this means that we have
$a_1\geq a_2 \geq \ldots \geq a_g$. The local system corresponding
to the standard representation defined by
$(a_1,\ldots,a_g,c)=(1,0,\ldots,0,1)$ is $R^1\pi_*{\QQ}(1)$.
If we do not specify $c$ then we assume that $c=\sum_{i=1}^g a_i$.
Duality now says that we have a non-degenerate pairing
$$
\VVl \times \VVl \to {\QQ}(-|\lambda|)
$$
with $|\lambda|=\sum \lambda_i$.

The irreducible representations of $M$ are parametrized by characters $\mu$
of $T$ which are $M$-dominant (i.e., with non-negative scalar
product with $\Phi_M^{+}$).
To such a character we can associate a locally free $O_{\Ag}$-module
(or vector bundle) $W_{\mu}$.
For example, the Hodge bundle of the universal family corresponds
to the representation $\gamma=(a,b,0,d) \mapsto \nu(\gamma)^{-1}a$ acting on
${\CC}^g$. Duality for $W_{\mu}$ says that $W_{\mu}^{\vee}=W_{-\sigma_1(w)}$
with $\sigma_1$ the longest element of $S_g$.
Faltings showed that one can extend the vector bundles
thus obtained to appropriate toroidal compactifications,
cf.\ \cite{F-C}, Thm.\ 4.2. 

\end{section}
\begin{section}{The BGG complex}\label{BGGcomplex}
Let $\Agf$ be a Faltings-Chai compactification of $\Ag$ and let $j: 
\Ag \to \Agf$
be the natural inclusion and let $i: D\to \Agf$ be the inclusion
of the divisor at infinity.
Recall that $D$ is a stack quotient of a Kuga-Satake variety,
namely a compactified quotient of the universal abelian variety of
dimension $g-1$ by the group ${\rm GL}(1,\ZZ)$ which acts 
by $\{\pm 1\}\in {\rm End}(X_{\eta})$  on the generic fibre $X_{\eta}$.

According to Deligne 
the logarithmic de Rham complex w.r.t.\ the divisor $D$ represents
$Rj_*{\CC}$. This generalizes for our sheaves $\VVl$, where the role of the
de Rham complex is played by the dual BGG complex which is obtained
by applying ideas of Bernstein-Gelfand-Gelfand of \cite{BGG} to our
situation as worked out in~\cite{F-C}. 

The dual BGG complex for $\lambda$ is a direct summand of the de Rham complex
for $\VVl^{\vee}$ and consists of a complex $K_{\lambda}^{\bullet}$
of vector bundles on $\Ag$:
$$
K_{\lambda}^{q} = \oplus_{w \in F_g, \ell(w)=q} W_{w * \lambda}^{\vee},
$$
where $w * \lambda = w(\lambda+\rho)-\rho$
and  $W_{\mu}^{\vee}=W_{-\sigma_1(\mu)}$. 
This complex is a filtered resolution of $\VVl^{\vee}$ on $\Ag$. 
The differentials
of this complex are given by homogeneous differential operators
on ${\DD}^{\vee}$.

The vector bundles $W_{\mu}$ extend over the compactification $\Agf$
and so do the differential operators, resulting in a complex
$\bar{K}_{\lambda}^{\bullet}$ on $\Agf$. We shall denote the extensions
of $W_{\mu}$ again by the same symbol $W_{\mu}$ (or by $\bar{W}_{\mu}$ 
if confusion might arise).
There is a variant of this complex $\bar{K}_{\lambda}^{\bullet}\otimes
O_{\Agf}(-D)$ and the differentials extend also for this complex.

The filtration on the dual BGG complex induce decreasing filtrations
on $\bar{K}_{\lambda}^{\bullet}$ and $\bar{K}_{\lambda}^{\bullet}\otimes O(-D)$
given by
$$
F^p(\bar{K}_{\lambda}^{\bullet})=
\oplus_{w \in F_g, f(w,\lambda)\geq p}\,  
\bar{W}_{w * \lambda}^{\vee}
$$
where $f(w,\lambda)= (\sum \lambda_i + \sum \mu_i)/2$ with
$\mu = -\sigma_1(w*\lambda)$,
and in an analogous way for $\bar{K}_{\lambda}^{\bullet}\otimes O(-D)$.
A term $W_{\mu}$ belongs to $F^p$ if and only if 
$(\sum \lambda_i +\sum \mu_i)/2 \geq p$.

In \cite{F-C}, p.\ 233, it is shown that the (filtered) dual BGG complex 
$\bar{K}_{\lambda}^{\bullet}$ is quasi-isomorphic to 
$Rj_*{\VV}_{\lambda}^{\vee}$,
while  $\bar{K}_{\lambda}^{\bullet}\otimes O(-D)$ is
quasi-isomorphic to $Rj_!\VV_{\lambda}^{\vee}$ and that the inclusion
$\bar{K}_{\lambda}^{\bullet}\otimes O(-D) \subset \bar{K}_{\lambda}^{\bullet}$
corresponds to the natural map $Rj_!\VVl^{\vee} \to Rj_*\VVl^{\vee}$.

We have an exact sequence of complexes
$$
0 \to \bar{K}_{\lambda}^{\bullet}\otimes O(-D) \to \bar{K}_{\lambda}^{\bullet}
\to \bar{K}_{\lambda}^{\bullet}|D \to 0.
$$
Therefore we can calculate the Eisenstein cohomology 
$e_{\rm Eis}(\Ag, \VVl^{\vee})$ by using the complex
$\bar{K}_{\lambda}^{\bullet}|D=i^* \bar{K}_{\lambda}^{\bullet}$ 
obtained by restricting the 
dual BGG complex
to $D$.
\end{section}

\begin{section}{A Calculation at the Boundary}\label{Leray}
The boundary stratum $D$ in $\Agf$ 
is a stratified space itself via the map $q$
to the Satake compactification. 
The part $D'=D\cap \Agp$ of $D$
lying in $\Agp$ has an \'etale cover
${\mathcal X}_{g-1}\to D'$
with ${\mathcal X}_{g-1} \to {\mathcal A}_{g-1}$
the universal principally polarized abelian variety of relative 
dimension $g-1$ and ${\mathcal X}_{g-1} \to {\mathcal A}_{g-1}$
factors through $q: D' \to {\mathcal A}_{g-1}/{\rm GL}(1,\ZZ)$.
The vector bundles $W_{\mu}$ extend canonically from $D'$ to $D$.

We shall calculate the rank $1$ part of the Eisenstein cohomology
by using the Leray spectral sequence for the complex 
$$
\oplus_{w \in F_g} W^{\vee}_{w*\lambda}
$$
by first restricting the factors $W^{\vee}_{w*\lambda}$ of the complex
$\oplus_{w \in F_g} W^{\vee}_{w*\lambda}$ to $D'$,
extending these to a suitable Faltings-Chai compactification $\bar{D}'$ and
then tensoring this complex with $O(-F)$ with $F$
the divisor at infinity $\bar{D}'-D'$ of $D'$ and by
calculating the cohomology using the Leray spectral sequence
for the map $q: D' \to {\mathcal A}_{g-1}/{\rm GL}(1,\ZZ)$.
(The tensoring with $O(-F)$ is done to get the rank-$1$ part 
of the Eisenstein cohomology and because of this the choice of
the compactification of $D'$ does not matter.)

As we shall now work at the same time with $\Ag$ and ${\mathcal A}_{g-1}$
we will write $W^{(g)}_{\mu}$ for $W_{\mu}$ on $\Agf$ in order
to avoid confusion. We first do a computation on ${\mathcal X}_{g-1}\to
{\mathcal A}_{g-1}$ and later take into account the action of
${\rm GL}(1,\ZZ)$.

\begin{proposition}\label{telescope}
For $W_{a}=W_{a}^{(g)}$ with $a=(a_1,\ldots,a_g)$ on $\tilde{\mathcal A}_g$ we have
$$
\sum_j (-1)^j R^jq_* (W_a^{(g)}|{\mathcal X}_{g-1})=
\sum_{k=1}^{g}
(-1)^{g-k} W^{(g-1)}_{(a_1,\ldots,a_{k-1},a_{k+1}-1,\ldots,a_{g}-1)}.
$$
\end{proposition}
\begin{proof}
We consider the Hodge bundle $\EE_g$ on $\Agf$. Its pullback to $\Xgmin$
fits into the exact sequence
$$
0 \to q^* \EE_{g-1} \to \EE_g \to O_{\Xgmin}\to 0,
$$
and we thus get
$$
Rq_* \EE_g = \EE_{g-1} \otimes Rq_*O_{\Xgmin} +Rq_*O_{\Xgmin}
= (\EE_{g-1}+1)\otimes \sum_{j=0}^{g-1} \wedge^j \EE_{g-1}^{\vee},
$$
since $R^iq_*O_{\Xgmin}=\wedge^i R^1q_*O_{\Xgmin}$.
Note that
$$
\sum_{j=0}^{g-1} \wedge^j \EE_{g-1}^{\vee}
=(-1)^k \sum_{k=0}^{g-1} W^{(g-1)}_{t_k}
$$
with $t_k$ denoting a vector $(0,\ldots,0,-1,\ldots,-1)$
of length $g-1$ with $g-1-k$ zeros. Since $W^{(g)}_a$ is made by
applying a Schur functor to the Hodge bundle
the exact sequence for $\EE_g$ implies
$$
Rq_* W_a=
{\rm Res}_{g-1}^g W^{(g)}_a \otimes \sum_{k=0}^{g-1} (-1)^k W^{(g-1)}_{t_k}.
$$
Here ${\rm Res}_{g-1}^g W_a$ is the bundle obtained from restriction
(branching) from ${\rm GL}(g)$ to ${\rm GL}(g-1)$. A well-known
formula from representation theory (cf.\ e.g.\ \cite{FH}) says that we thus get
$$
{\rm Res}_{g-1}^g W^{(g)}_a = \sum W^{(g-1)}_b,
$$
where the sum is over all (interlacing) $b=(b_1,\ldots,b_{g-1})$
with $a_1 \geq b_1 \geq a_2 \geq \cdots \geq b_{g-1}\geq a_g$.
Moreover, we have that $W^{(g-1)}(b) \otimes
\sum_{k=0}^{g-1} (-1)^k W^{(g-1)}_{t_k}$
is a signed sum of $W^{(g-1)}_{b'}$'s
with the sum running over all vectors $b'\in \ZZ^{g-1}$ obtained
from subtracting a $1$ from $k$ entries $b_i$ with $k$ between $0$ and $g-1$
and deleting those $b'$ that does not satisfy the condition that $b'_i \geq
b'_{i+1}$.
Carrying out the summation we see that most terms telescope away
and what remains is the right hand side of the statement
in the proposition. 
\end{proof}

For example, for $g=2$ we get that
$Rq_{*} W^{(3)}_{a,b,c}$ is equal to
$$
\sum_{\alpha,\beta}
W^{(2)}_{\alpha,\beta}-W^{(2)}_{\alpha,\beta-1}
-W^{(2)}_{\alpha-1,\beta}+W^{(2)}_{\alpha-1,\beta-1}
$$
where the sum is over all
$(\alpha,\beta)$ with
$a\geq \alpha \geq b \geq \beta \geq c$ and $W_2(r,s)=0$ if $r<s$. What remains
is $W_2(a,b)-W_2(a,c-1)+W_2(b-1,c-1)$.

\smallskip
However, we still need to take the action of ${\rm GL}(1,\ZZ)$ into account.
The non-trivial elements acts by $-1$ on the fibres of $q$. Therefore
only the terms in the right hand side of Proposition \ref{telescope} 
which are even (in the sense that $\sum_{i=1}^{k-1} a_i+\sum_{i=k+1}^g (a_i-1)
\equiv \, 0 (\bmod 2)$) will contribute to  
$$
\sum_j (-1)^j R^jq_* (W_a^{(g)}|D').
$$
\smallskip
Recall that we write $w * \lambda$ for the operation 
$\lambda \mapsto w(\lambda +\rho)-\rho$ of $W_g$ on the set of lambda's. 
Recall also we have $W(w*\lambda)^{\vee}=W(-\sigma_1 (w*\lambda))$ 
with $\sigma_1$ the longest element of $S_g$.

Our Eisenstein cohomology is now given by a complex which is
a sum over $w \in F_g$ of terms
$Rq_* W_{w * \lambda}^{\vee}$.
By Prop.\ \ref{telescope} the term $Rq_* W_{w * \lambda}^{\vee}$
yields a complex
$$
\sum_{l=1}^g (-1)^l W^{(g-1)}_{\tau_l(-\sigma(w*\lambda))},
$$
where $\tau_l$ applied to a vector $(a_1,\ldots, a_g)$ is the vector
$(a_1,\ldots,a_{l-1},a_{l+1}-1,\ldots,a_g-1)$.
We thus get a double sum $\sum_{w \in F_g} \sum_{l=1}^g X_{w,l}$ of terms
$X_{w,l}$ (that are (signed) sheaves of the form $W_{\mu}$)
which we rewrite as a sum of two double sums
$$
\sum_{k=1}^g \sum_{w\in F_g, w^{-1}(k)\leq k} X_{w,k} +
\sum_{k=1}^g \sum_{w\in F_g, w^{-1}(2g+1-k)\leq g} X_{w,k}.
$$
By Lemma \ref{subperms} the inner sum $\sum_{w\in F_g, w^{-1}(k)\leq k} X_{w,k}$
in the first double sum is equal to
$$
\sum_{u \in F_{g-1}} (-1)^{k} (W^{(g-1)}_{u*\tau^{\prime}_{k}(\lambda)})^{\vee},
$$
where $\tau_k^{\prime}(a_1,\ldots,a_g)=(a_1+1,\ldots,a_{k-1}+1,a_{k+1},
\ldots,a_g)$.
The double sum $
\sum_{k=1}^g \sum_{w\in F_g, w^{-1}(k)\leq k} X_{w,k}$
thus
contributes the complex
$$
K^{\prime}_{\lambda}:=\sum_{k=1}^g (-1)^k  K_{\tau^{\prime}_k(\lambda)}.
$$
The Hodge weight of these terms can be read off from Lemma 
\ref{subperms} and Faltings' results and thus contribute the cohomology of 
$\sum_k (-1)^k {\VV}_{\tau'_{k}(\lambda)}$.
(See the diagram in Section \ref{anotherexample} 
for an illustration in case $g=2$.)
In view of duality (Poincar\'e and Serre duality, see \cite{F-C}, p.\ 236)
the remaining $2^{g-1}g$  terms
of the sum
$\sum_{k=1}^g \sum_{w\in F_g, w^{-1}(2g+1-k)\leq g} X_{w,k}$
will contribute the dual terms:
$$
\sum_{k=1}^g (-1)^{k+1} K_{\tau^{\prime}_k(\lambda)}
\otimes \nu^{\lambda_k+g+1-k}
$$
and this contributes the cohomology of 
$\sum_k (-1)^{k+1} {\VV}_{\tau'_{k}(\lambda)}$ 
twisted by the power $\LL^{g+1+\lambda_k-k}$ of $\LL$.
But we need the rank $1$ part. To get this we consider the divisor
$\Delta$ at infinity of $\tilde{\mathcal A}_{g-1}$ 
(where $\tilde{\mathcal A}_{g-1}$ is defined as the closure
of the zero-section of ${\mathcal X}_{g-1}\to {\mathcal A}_{g-1}$)
and take 
$\bar{K}^{\prime}_{\lambda}\otimes O(-\Delta)$ instead of 
$K^{\prime}_{\lambda}$. Here $\bar{K}^{\prime}_{\lambda}$ is the extension
over $\tilde{\mathcal A}_{g-1}$ which we know to exist.
Then our Eisenstein cohomology is of the form
$$
\sum_{k=1}^g (-1)^k e_c({\mathcal A}_{g-1}, \VV_{\tau^{\prime}_k(\lambda)})
(1-\LL^{g+1+\lambda_k-k}).
$$
In the next section we show that this is compatible with the action
of the Hecke algebras.
\end{section}
%%%%%%%%%%%%%%%%
%%%%%%%%%%%%%%%%
\begin{section}{The Action of the Hecke Operators}
The Hecke algebra acts on the cohomology $H^*({\Ag},\VVl)$ and
$H^*_c(\Ag,\VVl)$ as explained in \cite{F-C}. It also acts in a 
compatible way on the dual BGG complex. These operators are
defined by algebraic cycles and this guarantees that they
respect the mixed Hodge structure on the Betti cohomology
and the Galois structure on \'etale cohomology. Moreover,
they are self-adjoint for Serre and Poincar\'e duality.

The compatible action of the Hecke operators on both 
$H^*({\Ag},\VVl)$ and $H^*_c(\Ag,\VVl)$ induces an action on
the (total) Eisenstein cohomology. We can see this action by means
of its action on the dual BGG complex and its restriction to
the boundary $D$. We now show that it factors through the action
of the Hecke algebra for ${\rm GSp}(2g-2,\ZZ)$.

The correspondences $T \to \Ag \times \Ag$ that define the Hecke 
operators extend in a natural way to $\Agp \times \Agp$. 
The (pullback of the) 
restriction of such a $T$ to ${\mathcal X}_{g-1}\times {\mathcal X}_{g-1}$,
a cover of $D'\times D'$,
is given by $T' \to {\mathcal X}_{g-1}\times {\mathcal X}_{g-1}$
which lies over a component of a Hecke correspondence $T^{\prime\prime}\to
{\mathcal A}_{g-1}\times {\mathcal A}_{g-1}$. In the next paragraph
we indicate this for the complex case.

The map
$T' \to T^{\prime\prime}$ is a universal family of abelian
varieties. The action is on bundles that are pullbacks from
${\mathcal A}_{g-1}$. We thus see that the action of $T'/T^{\prime \prime}$
is induced by an element of $\ZZ \subset {\rm End}(X_{\eta})$
with 
$X_{\eta}$, an abelian variety,
 the generic fibre of $T'$ over $T^{\prime \prime}$
on the cohomology of the stucture sheaf $O_{X_{\eta}}$.
This action is by scalars.
The action of $T^{\prime \prime}$ belongs to the action of the Hecke algebra
of ${\rm GSp}(2g-2,\ZZ)$. Hence the Hecke algebra of ${\rm GSp}(2g,\ZZ)$
on the Eisenstein cohomology factors through an action of the
Hecke algebra of ${\rm GSp}(2g-2,\ZZ)$. 

We work this out over $\CC$. A component of a Hecke correspondence
is defined by an embedding ${\mathcal H}_g \to 
{\mathcal H}_g \times {\mathcal H}_g$, with
${\mathcal H}_g$ the Siegel upper half plane, 
and given by equations
$$
wcz+dw-az-b=0, \eqno(1)
$$
where $w,z$ are in ${\mathcal H}_g$ and $(a,b;c,d)$ is an integral $2g\times 2g$-matrix
which lies in ${\rm GSp}(2g,\QQ)$.
This component can be extended to a correspondence for $\Agp$
that restricts to $T' \to  {\mathcal X}_{g-1}\times {\mathcal X}_{g-1}$
given by an embedding of ${\mathcal H}_{g-1} \times \CC^{g-1}$ into
$({\mathcal H}_{g-1} \times \CC^{g-1})^2$. We consider the behavior at infinity
given by 
$$
\lim_{t \to \infty} \left( \begin{matrix} z^{\prime} & \zeta \\
\zeta^t & it \\ \end{matrix} \right) \qquad z^{\prime} \in {\mathcal H}_{g-1},
\zeta \in \CC^{g-1}.
$$
In order that the component given by (1) does intersect our boundary
component the matrix $(a,b;c,d)$ must have the form
$$
\left(\begin{matrix} a' & 0 & b' & * \\
* & u  & * & * \\
c' & 0 & d' & * \\
0 & 0 & 0 & u^{-1} \\ \end{matrix}\right) 
\qquad (\left(\begin{matrix}
a' & b' \\ c' & d' \\ \end{matrix} \right) 
\in {\rm Sp}(2g-2,\QQ), u \in \QQ^*).
$$
This correspondence in ${\mathcal X}_{g-1}\times {\mathcal X}_{g-1}$
lies over a component of a Hecke correspondence $T^{\prime\prime}\to
{\mathcal A}_{g-1}\times {\mathcal A}_{g-1}$ given by
$$
w'c'z'+d'w'-a'z'-b'=0\qquad (z^{\prime}, w^{\prime} \in {\mathcal H}_{g-1}).
$$
and in the fibres it is given by
$$
\eta^t (c'z'+d')= u\zeta +(\alpha z'+\beta), 
$$
where $\eta$ is the analogue for $w$ of $\zeta$ for $z$.

\end{section}
%%%%%%%%%%%
%%%%%%%%%%%
\begin{section}{An Example: $g=1$}\label{g=1}
We consider the case of a local system $\VV_k={\rm Sym}^k(\VV)$
for $k$ even.
For $g=1$ the BGG complex is $0 \to j_*{\VV}_k^{\vee} \to W_{-k} \to W_{k+2}
\to 0$ with ${\VV}_k={\rm Sym}^k(R^1\pi_*(\QQ)(1))$. 
Similarly, there is a complex
$0 \to j_{!}{\VV}_k^{\vee} \to W_{-k}(-D) \to W_{k+2}(-D)
\to 0$ with $D$ the divisor $\tilde{\mathcal A}_1-{\mathcal A}_1$ at
infinity. By the exact sequence
$$
0 \to W_m(-D) \to W_m \to W_m|D \to 0
$$
we get
$$
e(j_*{\VV}_k^{\vee})-e(j_{!}{\VV}_k^{\vee})
=e(D,W_{-k}|D)-e(D,W_{k+2}|D).
$$
The bundle $W_{a}$ is associated to a representation of $Q$, the
parabolic subgroup. The action of the central multiplicative group
is not trivial: $W_{1}$ is associated to the representation where
$$
\left( \begin{matrix} a & b \\ 0 & d \\ \end{matrix} \right)
\text{ acts by multiplication by } 1/d \text{ on $\CC$}.
$$
In view of $\VV_{k}^{\vee}= \VV_k\otimes \nu^k$ we see that $e_{\rm Eis}
(\VV_k)$ is equals to $e(D,W_{-k}^{\prime}|D)-e(D,W_{k+2}^{\prime}|D)$,
where the prime refers to $k$ times twisting. This implies that
$W_{c}^{\prime}|D$ is $\CC((k+c)/2)$ on $D$ and that
$$
e_{\rm Eis}({\VV}_k)=\LL^0-\LL^{k+1}
$$
for $k\geq 0$ even.

Our answer for the Eisenstein cohomology
is compatible with Poincar\'e duality. It has a part with Hodge
weight $k+1$ and one with Hodge weight~$0$. The part of Hodge weight $k+1$
occurs for $k\geq 2$ in the exact sequence
$$
0 \to H^0(W_{k+2}\otimes O(-D))
\to H^0(W_{k+2}) \to H^0(D,W_{k+2}{|D})
\to 0.
$$
This can be identified with
$$
0 \to S_{k+2} \to M_{k+2} \to \CC(k+1) \to 0,
$$
with $M_{k+2}$ (resp.\ $S_{k+2}$) the space of
modular forms (resp.\ of cusp forms) of weight $k+2$ on ${\rm SL}(2,\ZZ)$.
The other part occurs in
$$
0 \to H^0(D,W_{-k}{|D})
\to H^1(W_{-k}\otimes O(-D))
\to H^1(W_{-k}) \to 0.
$$
All in all we get
$$
e_c({\mathcal A}_1,\VV_{k})
= -S_{k+2} \oplus \bar{S}_{k+2}-1
$$
and
$$
e({\mathcal A}_1,\VV_{k})=-S_{k+2} \oplus \bar{S}_{k+2}-\CC(k+1).
$$
We identify $\CC$ and $\CC(k+1)$ with the Eisenstein cohomology. Note that if
we let $W_m^{\prime}=W_m\otimes O(-\Delta)$ with $\Delta$ the divisor on $D$
that is the fibre over the cusp $\infty$ of $\tilde{\mathcal A}_1$
we also have the identities
$$
\begin{aligned}
e_c({\mathcal A}_1,\VV_k)&= 
-e(\tilde{\mathcal A}_1,W_{k+2}^{\prime})+e(\tilde{\mathcal A}_1,W_{-k}^{\prime})
\cr
e({\mathcal A}_1,\VV_k)&= 
-e(\tilde{\mathcal A}_1,W_{k+2})+e(\tilde{\mathcal A}_1,W_{-k}).
\cr
\end{aligned}
$$

For even $k\geq 4$ we let $S[k]$ denote the motive of cusp forms of weight $k$
on ${\rm SL}(2,\ZZ)$ as constructed by Scholl, cf.\ \cite{Scholl}, see also
\cite{C-F}. For $k=2$ we put $S[2]=-\LL-1$. In the category of Hodge structures
we have for $k\geq 2$
$$
S[k+2]= e(\tilde{\mathcal A}_1,W^{\prime}_{k+2}) -e(\tilde{\mathcal A}_1,W_{-k}).
$$
Note that by Serre duality we have $H^1(\tilde{\mathcal A}_1,W_{-k})\cong
H^0(\tilde{\mathcal A}_1,\Omega^1\otimes W_k)^{\vee}=H^0(\tilde{\mathcal A}_1,
W_{k+2}^{\prime})^{\vee}$ but if we take into account the action of the central
$\GG_m$ we have to twist by $\eta^{k+1}$.
We now have for even $k\geq 0$
$$
e_c({\mathcal A}_1,\VV_{k})= -S[k+2]-1, \qquad e({\mathcal A}_1,\VV_{k})
=-S[k+2]-\LL^{k+1}.
$$
\end{section}
%%%
%%%
\begin{section}{Another Example: $g=2$}\label{anotherexample}
We now look at the case $g=2$ and consider the local system $\VV_{l,m}$
with $l\equiv m \, (\bmod 2)$. Calculations with Carel Faber in \cite{FvdG}
led to the formulas for this case.
We have a standard Faltings-Chai compactification
$\tilde{\mathcal A}_2$ in this case which coincides with Igusa's blow-up of the Satake compactification and also with the moduli space $\overline{\mathcal M}_2$
of stable curves of genus $2$. 
We consider the full Eisenstein cohomology
$$
e_{\rm Eis}({\mathcal A}_2,\VV_{l,m}):=
e({\mathcal A}_2,Rj_{*}\VV_{l,m})-e({\mathcal A}_2,Rj_{!}\VV_{l,m}).
$$

\begin{theorem}
The Eisenstein cohomology $e_{\rm Eis}({\mathcal A}_2,\VV_{l,m})$
is given by
$$
\begin{aligned}
-s_{l-m+2}&(1-\LL^{l+m+3})+  s_{l+m+4}(\LL^{m+1}-\LL^{l+2})+ \qquad \\
\qquad + & \begin{cases}
e_c({\mathcal A}_1,(\VV_m)(1-\LL^{l+2})-(\LL^{l+2}-\LL^{l+m+3})
   & \text{ $l$ even,} \cr
-e_c({\mathcal A}_1,\VV_{l+1})(1-\LL^{m+1})-(1-\LL^{m+1}) & \text{$l$ odd.} \cr
\end{cases} \\
\end{aligned}
$$
\end{theorem}
\smallskip
\noindent
Alternatively, the Eisenstein cohomology can be written as
$$
\begin{aligned}
-(s_{l-m+2}+1)(1-\LL^{l+m+3})+ & s_{l+m+4}(\LL^{m+1}-\LL^{l+2})+ \qquad \\
\qquad + & \begin{cases}
-S[m+2](1-\LL^{l+2}) & \text{ $l$ even,} \cr
S[l+3] (1-\LL^{m+1}) & \text{ $l$ odd.} \cr
\end{cases} \\
\end{aligned}
$$
We can check this for example for $l=m=0$. We have
$$
e(Rj_{*}\VV_{0,0})-e(Rj_{!}\VV_{0,0}) = 1 +\LL-\LL^2-\LL^3,
$$
while for the {\sl compactly supported Eisenstein cohomology} (i.e., the kernel
of $H^*_c \to H^*$) we find $1+\LL$, which fits.
If $\VV_{l,m}$ is a regular local system (i.e., $l>m>0$) then the Hodge weights
of the terms are either $>l+m+3$ or $<l+m+3$ and this determines whether 
the term belongs to compactly supported Eisenstein cohomology, cf.,
\cite{P}{ Thm.\ 3.5}. We thus get the result announced in \cite{FvdG}.

\begin{corollary}\label{gistwee}
(\cite{FvdG})
The compactly supported Eisenstein cohomology for a regular local
system $\VV_{l,m}$ is given by
$$
s_{l-m+2} -s_{l+m+4} \LL^{m+1} + \begin{cases}
S[m+2]+1 & \text{ $l$ even,} \cr
-S[l+3] & \text{$l$ odd.} \cr
\end{cases}
$$
\end{corollary}
In \cite{FvdG} one finds numerical confirmation of these formulas.
The BGG complex for $j_*\VV_{l,m}^{\vee}$ is:
$$
0 \to j_*\VV_{l,m}^{\vee} \to W_{-m,-l}\to W_{m+2,-l}
\to W_{l+3,1-m} \to W_{l+3,m+3} \to 0
$$
and for the compactly supported cohomology
the similar complex is
$$
0 \to j_!\VV_{l,m}^{\vee} \to W_{-m,-l}\otimes O(-D) \to \ldots
$$
The extended complexes over $\tilde{\mathcal A}_g$ are quasi-isomorphic to
$Rj_*(\VVl)^{\vee}$ and $Rj_!(\VVl)^{\vee}$, see \cite{F-C}{ Prop.\ 5.4}
and above.
We have to twist $l+m$ times if we work with $Rj_*(\VVl)$ and
$Rj_!(\VVl)$. We hope that the details of the first part of the
proof in this case will illustrate the proof of the formula for
$e_{\rm Eis,1}$ for the general case. The short exact sequence
$$
0 \to W_{\mu}(-D) \to W_{\mu} \to W_{\mu}|D \to 0
$$
yields
$$
\begin{aligned}
e(Rj_*\VV_{l,m})-e(j_!\VV_{l,m})&=
e(W_{-m,-l}|D)-e(W_{m+2,-l}|D)+\\
& +e(W_{l+3,1-m}|D)-e(W_{l+3,m+3}|D))\\
\end{aligned}
$$
if we take into account a $l+m$ th twist.
Therefore we first determine the Euler characteristic $e(W_{a,b}|D)$.
Note that $D$ is the quotient by $-1$ of the compactification 
$\tilde{\mathcal X}_1 \to \tilde{\mathcal A}_1$ of
the universal elliptic curve ${\mathcal X}_1 \to {\mathcal A}_1$.
We stratify $\tilde{\mathcal X}_1$ by the open part and the fibre $F$
over the cusp $\infty$ of $\tilde{\mathcal A}_1$. The cohomology
$$
e(\tilde{\mathcal A}_{1},\sum_j (-1)^j R^jq_{*} (W_{\mu}|{\mathcal X}_{1}))
$$
can be calculated via the exact sequence for the Hodge bundle $\EE_2=W_{1,0}$
$$
0 \to q^*\EE_{1} \to \EE_2|{\mathcal X}_{1} \to O_{{\mathcal X}_{1}}
\to 0
$$
More generally, the pull back of $W_{a,b}$
to ${\mathcal X}_1$ is ${\rm Sym}^{a-b} \EE_1 \otimes \det \EE_1^b$
with $\EE_1$ the Hodge bundle. But we need to keep track of the twisting.
The exact sequence above gives (using that $R^iq_*O_{{\mathcal X}_{1}}=
\wedge^i \EE_{1}^{\vee}$)
$$
Rq_*W_{a,b}|{\mathcal X}_1 = \sum_{a\geq \nu \geq b} W_{\nu} 
\otimes (1-\EE_1^{\vee})
$$
and carrying out the summation we find
$$
Rq_*W_{a,b}|{\mathcal X}_1 = W_{a} -W_{b-1}.
$$
We now replace the $W_{\mu}$'s in our calculation by 
$W_{\mu}':=W_{\mu}(-F)$ and then calculate separately the contribution
from the stratum $F$.
We  collect the terms. An element $w\in F_2\subset S_4$ is identified
with its images $[w(1) w(2)]$ of the elements $1$ and $2$.
%%A contribution
%%$e(\tilde{\mathcal A}_1,W_{\mu}$ coming from $w \in F_g$ should 
%%have Hodge weight $0$ (resp.\ $m+1,l+2,l+m+3$) for $\ell(w)=0$
%%resp.\ $1,2,3$). Since the Hodge weight of 
%%$e(\tilde{\mathcal A}_1,W_k^{\prime})$ 
%%is $0$ for $k<-2$ and $k+1$ for $k>2$ we see with which power of
%%$\CC(1)$ we must twist.

\begin{center}
\begin{tabular}[t]{|c|c|c|c|c|c|}
\hline
$w$ & $\ell(w)$  & $-\sigma_1(w*\lambda)$ & contribution \\
\hline
$[12]$ & $0$ & $(-m,-l)$ & {\color{red}$W_{-m}^{\prime}$}
{\color{blue}$-W_{-l-1}^{\prime}$} \\
$[13]$ & $1$ & $(m+2,-l)$ & {\color{red}$-W_{m+2}^{\prime}$}
$+W_{-l-1}^{\prime}(\nu^{m+1})$\\
$[24]$ & $2$ & $(l+3,1-m)$ & {\color{blue}$W_{l+3}^{\prime}$}
{\color{green}$-W_{-m}^{\prime}(\nu^{l+2})$}\\
$[34]$ & $3$ & $(l+3,m+3)$ & $-W_{l+3}^{\prime}(\nu^{m+1})$
{\color{green}$+W_{m+2}^{\prime}(\nu^{l+2})$}\\
\hline
\end{tabular}
\end{center}

\bigskip
\noindent
We can collect this into $-e_c({\mathcal A}_1,\VV_m)(1-\LL^{l+2})+
e_c({\mathcal A}_1,\VV_{l+1})(1-\LL^{m+1})$.

\smallskip
We now treat the codimension $2$ boundary contribution.
It is a priori clear that the outcome will be a polynomial in $\LL$
as we are calculating over a toric curve;
in fact one can also deduce that the exponents of $\LL$
that occur are in $\{l+m+3,l+2,m+1,0\}$.

Let $F$ be the fibre of $q$
over $\tilde{\mathcal A}_0$. Note that $F$ has dimension~$1$.
The neighborhood of $F$ in the toroidal compactification
$\tilde{\mathcal A}_2$ is a toric variety obtained by glueing
infinitely many copies of affine $3$-space $A^3$ and dividing
through an action of ${\rm GL}(2,\ZZ)$;
more precisely it is an orbifold constructed as follows. The cone of
symmetric real positive definite $2\times 2$ matrices has a natural
cone decomposition invariant under the action of ${\rm GL}(2,\ZZ)$.
The group ${\rm GL}(2,\ZZ)$ acts on the cone of positive definite
$2 \times 2$ matrices by
$$
C=
\left( \begin{matrix} \alpha & \beta \\ \beta & \gamma \\
\end{matrix} \right)
\mapsto A^tCA
$$
The cone decomposition consists of the orbit under ${\rm GL}(2,\ZZ)$
of the cone spanned by the three matrices
$$
\left( \begin{matrix} 1 & 0 \\ 0 & 0 \\ \end{matrix} \right), \quad
\left( \begin{matrix} 0 & 0 \\ 0 & 1 \\ \end{matrix} \right), \quad
\left( \begin{matrix} 1 & 1 \\ 1 & 1 \\ \end{matrix} \right).
$$
Associated to this we have the toroidal variety ${\mathcal T}$
with ${\rm GL}(2,\ZZ)$-action.
The coordinate axes in the $A^3$'s are glued to a union of $\PP^1$'s.
The $1$-skeleton ${\mathcal T}^1$ of the quotient of this toric variety
under ${\rm GL}(2,\ZZ)$ can be identified with $F$.
In general in level $n\geq 3$ we have $(1/4)n^3\prod_{p|n} (1-p^{-2})$
copies of $\PP^1$ meeting three at each one of the
$(1/6)n^3 \prod_{p|n}(1-p^{-2})$ points.
For $n=3$ this looks like a tetrahedron, a cube for $n=4$ etc.\
and in general as the polyhedral
decomposition of the Riemann surface $\Gamma(n)\backslash {\mathcal H}_1$
with $\Gamma(n)$ the full level $n$ congruence subgroup.

The associated quadratic form $\alpha x^2+2\beta xy + \gamma y^2$
determines a point $ z=(-\beta + \sqrt{\beta^2-\alpha \gamma})/\alpha$
in the upper half plane (this factors through scaling by positive reals).
The action of ${\rm diag}( 1,-1)$ is given by $\beta \mapsto -\beta$
and induces $z\mapsto -\bar{z}$ on the upper half plane.

In order to calculate the complex $\bar{K}_{\lambda}|F$ 
we restrict the Hodge bundle $\EE$ to $F$. Note that the fundamental domain
for the action of $\Gamma[n]\subset {\rm GL}(2,\ZZ)$
is a cone over a fundamental domain for the modular curve of level $n$.

The restriction of the Hodge bundle to a $\PP^1$ in $F$
is of the form $O_{\PP^1}\oplus O_{\PP^1}$. The pullback of $\EE$ to
${\mathcal T}^1$ is flat vector bundle of rank $2$ determined by
the standard representation of ${\rm GL}(2)$. So the bundle
$W(a,b)$ is associated to the irreducible representation
of type $(a,b)$ of ${\rm GL}(2)$. This and the toroidal
construction makes it possible to express the cohomology in terms of
group cohomology.
\begin{proposition}  Let
$$
U=W_{-m,-l}-W_{m+2,-l}+W_{l+3,1-m}-W_{l+3,m+3}.
$$
Then the Euler characteristic $e_c(F,U|F)$
equals
$$
-s_{l-m+2}(\LL^0 -\LL^{l+m+3})+ s_{l+m+4}(\LL^{m+1}-\LL^{l+2})+
\begin{cases}
-\LL^{l+2}+\LL^{l+m+3} &  \text{ $l$ even}  \\
-1+\LL^{m+1} &  \text{ $l$ odd.} \\
\end{cases}
$$
\end{proposition}

\noindent
\begin{proof}
The Euler characteristic $e_c(F,W(a,b))$ can be expressed
in group cohomology. We calculate it by using the stratification of $F$
by the open part and the cusp. We then get the Euler characteristic
of compactly supported cohomology of ${\rm GL}(2,\ZZ)$ with values in $V_{a-b}$.
This gives for the four cases $(-m,-l), \ldots, (l+3,m+3)$
the contributions $(-s_{l-m+2}-1)\LL^x$, $(-s_{l+m+4}-1)\LL^y$,
$(s_{l+m+4}+1)\LL^{l+m+3-y}$ and
$(s_{l-m+2}+1)\LL^{l+m+3-x}$ for appropriate $x$ and $y$
which turn out to be $0$ and $m+1$.
We have to add the contribution from
the cusp. This gives
$\LL^{l+2}-\LL^{l+m+3}$ for $l$ odd
and $1-\LL^{m+1}$ for $l$ even. 
For this just look at the action of ${\rm diag}(1,-1)$. 
If $\ell$ (and hence $m$) is even only the first two terms
($W_{-m,-l}$ and $-W_{m+2,-l}$)
contribute because of the sign, while for the odd case the
two other terms contribute.
\end{proof}
\begin{remark}
The same method suffices to prove the analogues for the moduli
space ${\mathcal A}_2[n]$ of abelian surfaces with a level $n$
structure. 
\end{remark}

\end{section}
%%%

\end{document}